\documentclass[12pt]{article}

\usepackage{amsmath}
\usepackage{amssymb}
\usepackage{color}
\usepackage[dvips]{graphics}
\usepackage{graphicx}
\usepackage{cite}

\date{}
\title{New representations of Pad\'{e} and Pad\'{e}--type approximants}
\author{Claude Brezinski\thanks{Laboratoire Paul Painlev\'e, UMR CNRS 8524, UFR
de Math\'ematiques, Universit\'e des Sciences
et Technologies de Lille,
59655--Villeneuve d'Ascq cedex,
France, E--mail: {\tt Claude.Brezinski@univ-lille1.fr}.}
\and Michela Redivo--Zaglia\thanks{Universit\`a degli Studi di Padova,
Dipartimento di Matematica,
Via Trieste 63, 35121--Padova,
Italy. E--mail: {\tt Michela.RedivoZaglia@unipd.it}.
}
}

\begin{document}

\maketitle

\vskip 1mm

\centerline{{\it In memoriam nostri}  Pablo Gonz\'alez Vera}

\begin{abstract}
Pad\'{e} approximants are rational functions whose series expansion match a given series as far as possible. These approximants are usually written under a rational form. In this paper, we will show how to write them also under two different barycentric forms, and under a partial fraction form, depending on free parameters. According to the choice of these parameters, Pad\'{e}--type approximants can be obtained under a barycentric or a partial fraction form.
\end{abstract}

\noindent {\bf Keywords:} Pad\'e approximation, barycentric rational function, partial fraction.

\section{Introduction}

This paper describes new mathematical expressions for Pad\'e approximants, and some of their variants.
A Pad\'{e} approximant is a rational function whose power series expansion in ascending powers of the variable matches a given formal power series as far as possible \cite{bak,birk}. Thus, it can be understood as a rational Hermite interpolant at zero, and it is usually written under the form of a rational fraction or as the convergent of a certain continued fraction.
On the other hand, a rational interpolant can be given under the form of a rational fraction, or as the convergent of a continued fraction, or under a barycentric form.

\vskip 2mm

In this paper, we will show that a Pad\'{e} approximant can also be written under (at least) two different barycentric rational forms which depend on arbitrary parameters. Such an approximant will be called a {\it barycentric Pad\'{e} approximant} (in short {\sc bpa}).
According to the choice of these free parameters, Pad\'{e}--type approximants are also obtained under this form and we call them {\it barycentric Pad\'{e}--type approximant} (in short {\sc bpta}).
Then, we will show how to write a Pad\'e approximant under a partial fraction form,
called a {\it partial fraction Pad\'e approximant} (in short {\sc pfpa}).
The case of partial Pad\'e approximants \cite{ppa}, where some poles and/or zeros are imposed, could be treated similarly.

\section{Rational form}

Let $f$ be a formal power series
\begin{equation}
\label{f}
f(t)=c_0+c_1t+c_2t^2+\cdots
\end{equation}
We consider the rational function
\begin{equation}
\label{r1}
R_{p,q}(t)=\frac{\displaystyle \sum_{i=0}^p a_it^i}{\displaystyle \sum_{i=0}^q b_it^i}.
\end{equation}
If the $b_i$'s are arbitrarily chosen (with $b_0b_q \neq 0$), and if the $a_i$'s are computed by
\begin{equation}
\left.
\begin{array}{rcl}
a_0 & = & c_0b_0 \\
a_1 & = & c_1b_0 + c_0b_1 \\
&\vdots &  \\
a_p & = & c_pb_0 + c_{p-1}b_1 + \cdots + c_{p-q}b_q
\end{array}
\right\}
\label{ai}
\end{equation}
with the convention that $c_i = 0$ for $i < 0$, then $R_{p,q}$ is the {\it Pad\'{e}--type approximant} of $f$ \cite{pta}, it is denoted by $(p/q)_f$, and it holds
$$(p/q)_f(t)-f(t)={\cal O}(t^{p+1}).$$
This {\it accuracy--through--order condition} means that the first $p+1$ coefficients of the power series expansion of $(p/q)_f$ in ascending powers of the variable $t$ match those of the series $f$.

\vskip 2mm

Moreover, if the $b_i$'s are taken as the solution of the system
\begin{equation}
\left.
\begin{array}{rcl}
0&=& c_{p+1}b_0+ c_{p}b_1 + \cdots + c_{p-q+1}b_q \\
&\vdots &  \\
0&=& c_{p+q}b_0+ c_{p+q-1}b_1 + \cdots + c_{p}b_q,
\end{array}
\right\}
\label{bi}
\end{equation}
with $b_0=1$ (a rational function is defined up to a multiplying factor), then $R_{p,q}$ is the {\it Pad\'{e} approximant} of $f$
\cite{bak,birk}, it is denoted by $[p/q]_f$, and it holds
$$[p/q]_f(t)-f(t)={\cal O}(t^{p+q+1}).$$
Thus, the first $p+q+1$ coefficients of the series expansion of $[p/q]_f$ are identical to those of $f$.
Moreover, we have
$$
[p/q]_f(t)=\left|
\begin{array}{cccc}
t^qf_{p-q}(t) & t^{q-1}f_{p-q+1}(t) & \cdots & f_{p}(t) \\
c_{p-q+1} & c_{p-q+2} & \cdots & c_{p+1} \\
\vdots & \vdots & & \vdots \\
c_{p} & c_{p+1} & \cdots & c_{p+q}
\end{array}
\right| \Big/ \left|
\begin{array}{cccc}
t^q & t^{q-1} & \cdots & 1 \\
c_{p-q+1} & c_{p-q+2} & \cdots & c_{p+1} \\
\vdots & \vdots & & \vdots \\
c_{p} & c_{p+1} & \cdots & c_{p+q}
\end{array}
\right|,
$$
with $f_n(t)=c_0+\cdots+c_nt^n$, the $n$th partial sum of the series $f$ ($f_n$ is identically zero for $n<0$).

\vskip 2mm

In the case of a {\it partial Pad\'{e} approximant}, a part of its numerator and/or its denominator is arbitrarily chosen, and the remaining part is taken so that its power series expansion matches $f$ as far as possible \cite{ppa}.

\section{Barycentric forms ({\sc bpa})}

In this section, we consider rational functions written under two different barycentric forms

\begin{equation}
\label{f1}
\mbox{Form 1:} \qquad \; \widetilde R_{p,q}(t)=\frac{\displaystyle \sum_{i=0}^p \frac{\widetilde a_i}{\widetilde p_i-t}}{\displaystyle \sum_{i=0}^q \frac{\widetilde b_i}{\widetilde z_i-t}}
\end{equation}
or
\begin{equation}
\label{f2}
\mbox{Form 2:} \qquad \widehat R_{p,q}(t)=\frac{\displaystyle \sum_{i=0}^p \frac{\widehat a_i}{1-\widehat p_it}}{\displaystyle \sum_{i=0}^q \frac{\widehat b_i}{1-\widehat z_it}},
\end{equation}
where the $\widetilde p_i$'s, the $\widetilde z_i$'s, the $\widehat p_i$'s, and the $\widehat z_i$'s are given points in the complex plane. We assume that all the $\widetilde p_i$'s are distinct, and also all the $\widetilde z_i$'s,
all the $\widehat p_i$'s, and all the $\widehat z_i$'s.

Obviously, the forms \eqref{f1} and \eqref{f2} can be deduced one from each other by setting
$\widetilde a_i=\widehat a_i/\widehat p_i$, $\widetilde b_i=\widehat b_i/\widehat z_i$, $\widetilde p_i=1/\widehat p_i$,
and $\widetilde z_i=1/\widehat z_i$, for all points different from zero. In the sequel, when it is not necessary to distinguish between the two forms and when it is possible to treat them simultaneously, any of them will be simply denoted by $R_{p,q}$, and the parameters by $a_i$, $b_i$, $p_i$, and $z_i$ respectively.

\vskip 2mm

In both cases, we want to determine the coefficients $a_i$ and $b_i$ such that
\begin{equation}
\label{cond}
R_{p,q}(t)-f(t)={\cal O}(t^{p+q+1}).
\end{equation}
Due to this property, and although $R_{p,q}$ is not always identical to the Pad\'e
approximant $[p/q]_f$ of the series $f$ as we will see below,
such a rational function will be called a {\it barycentric Pad\'{e} approximant} and denoted {\sc bpa}.

\vskip 2mm

Before explaining how to compute the coefficients of such an approximant, let us begin by some important remarks:

\begin{enumerate}

\item
It is easy to see that, for \eqref{f1}, the $\widetilde p_i$'s are poles of $\widetilde R_{p,q}$ and the $\widetilde z_i$'s are zeros of it while, for \eqref{f2}, it is the $1/\widehat p_i$'s and the $1/\widehat z_i$'s which play these roles. Therefore, if some poles and zeros of $f$ are known, they can be introduced into the construction of the approximant as in the case of partial Pad\'{e} approximants \cite{ppa}. For the form \eqref{f1}, we will assume that $\forall i, \widetilde p_i \neq 0, \widetilde z_i \neq 0$.
The reason for this condition will be clearly seen in Section \ref{s311}. If, in \eqref{f2}, some $\widehat p_i$'s and/or some
$\widehat z_i$'s are zero the degree of the numerator and/or the degree of the denominator reduces accordingly.

\item
After reducing the sum in the numerator of \eqref{f1} to its common denominator and also the sum in the denominator, $\widetilde R_{p,q}$ becomes
$$\widetilde R_{p,q}(t)=\frac{N_p(t)\prod_{i=0}^q(\widetilde z_i-t)}{D_q(t)\prod_{i=0}^p(\widetilde p_i-t)},$$
where $N_p$ is a polynomial of degree $p$ and $D_q$ a polynomial of degree $q$.

If $\forall i, \widetilde p_i \neq \widetilde z_i$, then $\widetilde R_{p,q}$ has a numerator and a denominator both of degree
$p+q+1$  at most. Thus, the order of approximation of $\widetilde R_{p,q}$ is one less than the order of approximation of the Pad\'e--type approximant with the same degrees \cite{pta}. We will discuss below how to improve this order.
The $\widetilde p_i$'s and the $\widetilde z_i$'s can be selected so that $\widetilde R_{p,q}$ possesses other interesting properties such as, for example, the preservation of as many moments of $f$ as possible. However, they cannot be chosen after the $\widetilde a_i$'s and the $\widetilde b_i$'s have been computed since, as we will see below, these coefficients depend on them.

If some of the $\widetilde p_i$'s coincide with some of the $\widetilde z_i$'s,   then a cancelation occurs and it lowers the degrees accordingly. If, when $p<q$, $\widetilde p_i=\widetilde z_i$ for $i=0,\ldots,p$, the product in the denominator disappears and the product in the numerator reduces to $\prod_{i=p+1}^q(\widetilde z_i-t)$. Thus $\widetilde z_{p+1},\ldots,\widetilde z_q$ are zeros of $\widetilde R_{p,q}$.
When $q<p$ and $\widetilde p_i=\widetilde z_i$ for $i=0,\ldots,q$, it is the product in the numerator which disappears and the product in the denominator reduces to $\prod_{i=q+1}^p(\widetilde p_i-t)$. Thus $\widetilde p_{q+1},\ldots,\widetilde p_p$ are poles of $\widetilde R_{p,q}$.

Similar remarks hold for \eqref{f2}.

\item
In both cases, if $\forall i, p_i = z_i$ and $p=q$, then $R_{p,p}$
has a numerator and a denominator both of degree $p$ at most. Thus, thanks to the condition \eqref{cond},
$R_{p,p}$ is the usual Pad\'{e} approximant $[p/p]_f$ of $f$, and it is such that $R_{p,p}(t)-f(t)={\cal O}(t^{2p+1})$. Due to its uniqueness, this approximant is, in theory, independent of the choice of the $p_i$'s. However, in practice, the choice of the $p_i$'s can influence the stability of the approximant, an important issue yet to be studied.

Barycentric Pad\'e approximants (which, in this case, are true Pad\'e approximants) with arbitrary degrees in the numerator and in the denominator can be constructed as follows.
Let us write $f$ as
$$f(t)=c_0+\cdots+c_{n-1}t^{n-1}+t^nf^n(t) \quad \mbox{with} \quad f^n(t)=c_n+c_{n+1}t+\cdots$$
The approximant
$$R_{n+p,p}(t)=c_0+\cdots+c_{n-1}t^{n-1}+t^nR_{p,p}(t),$$
where $R_{p,p}(t)$ is now the barycentric Pad\'e approximant of the series $f^n$, satisfies
$$R_{n+p,p}(t)-f(t)={\cal O}(t^{n+2p+1}).$$
Thus it is identical to the Pad\'e approximant $[n+p/p]_f$ independently of the choice of the $p_i$'s.
Similarly, write $f$ as
$$f(t)=t^{-n}f^{-n}(t) \quad \mbox{with} \quad f^{-n}(t)=0+0t+\cdots+0t^{n-1}+c_0t^n+c_1t^{n+1}+\cdots.$$
The approximant
$$R_{p,p+n}(t)=t^{-n}R_{p+n,p+n}(t),$$
where $R_{p+n,p+n}$ is the barycentric Pad\'e approximant of the series $f^{-n}$ satisfies
$$R_{p,n+p}(t)-f(t)={\cal O}(t^{n+2p+1}).$$
Thus it is identical to the Pad\'e approximant $[p/n+p]_f$ independently of the choice of the $p_i$'s.

\item
Let us remind that if, in \eqref{f1}, $p=q$, $\forall i, \widetilde p_i=\widetilde z_i$ and $\widetilde a_i=w_i f(\widetilde p_i), \widetilde b_i=w_i \neq 0$, then $\widetilde R_{p,p}(p_i)=f(\widetilde p_i)$ independently of the choice of the $w_i$ \cite{csww}. Thus, the $w_i$'s can be chosen so that, in addition, $\widetilde R_{p,p}$
matches the series $f$ as far as possible as proposed in \cite{ptri}.
\end{enumerate}

\subsection{Computation of the coefficients}
\label{cc}

Since the accuracy--through--order condition \eqref{cond} contains $p+q+1$ relations while $R_{p,q}$ has $p+q+2$ coefficients to be determined, an additional condition needs to be imposed. It is the so--called {\it normalization condition}. Because $R_{p,q}$ approximates $f$ around zero, its denominator should not vanish at this point. Thus, since a rational function is determined apart a common multiplying factor in its numerator and its denominator, it is convenient to choose the normalization condition for \eqref{f1} and \eqref{f2} as
\begin{equation}
\label{norma}
\sum_{i=0}^q \frac{\widetilde b_i}{\widetilde z_i}=1 \mbox{~~for~~}  \eqref{f1},\quad \mbox{and} \quad \sum_{i=0}^q \widehat b_i=1  \mbox{~~for~~}  \eqref{f2}.
\end{equation}

\subsubsection{Form 1}
\label{s311}

For the form \eqref{f1}, having chosen the $\widetilde p_i$'s and the $\widetilde z_i$'s, the accuracy--through--order condition \eqref{cond} can be written
$$\sum_{i=0}^p \frac{\widetilde a_i/\widetilde p_i}{1-t/\widetilde p_i}=(c_0+c_1t+c_2t^2+\cdots)\sum_{i=0}^q\frac{\widetilde b_i/\widetilde z_i}{1-t/\widetilde z_i}.$$
But, we have $1/(1-t/\widetilde p_i)=1+t/\widetilde p_i+t^2/\widetilde p_i^2+\cdots$, and a similar expansion for $1/(1-t/\widetilde z_i)$. Thus, after replacement, the preceding relation becomes
$$\sum_{i=0}^p \frac{\widetilde a_i}{\widetilde p_i}\left(1+\frac{t}{\widetilde p_i}+\frac{t^2}{\widetilde p_i^2}+\cdots\right)=(c_0+c_1t+c_2t^2+\cdots)
\sum_{i=0}^q \frac{\widetilde b_i}{\widetilde z_i}\left(1+\frac{t}{\widetilde z_i}+\frac{t^2}{\widetilde z_i^2}+\cdots\right).$$

Identifying the coefficients of the identical powers of $t$ in both  sides, and taking into account the normalization condition leads to the system of equations for computing the coefficients ${\widetilde a_i}$ and ${\widetilde b_i}$
\begin{equation}
\label{coef1}
\left.
\begin{array}{l}
\displaystyle \sum_{i=0}^q \frac{\widetilde b_i}{\widetilde z_i}=1 \\
\displaystyle \sum_{i=0}^p \frac{\widetilde a_i}{\widetilde p_i^{k+1}}-\sum_{i=0}^q \widetilde b_i \sum_{j=0}^{k-1} \frac{c_j}{\widetilde z_i^{k-j+1}}=c_k, \qquad k=0,\ldots,p+q.
\end{array}
\right\}
\end{equation}

Obviously, the sum on $j$ is empty for $k=0$.

\vskip 2mm

Thus, the coefficients of the series expansion of $\widetilde R_{p,q}(t)=\widetilde d_0+\widetilde d_1t+\widetilde d_2t^2+\cdots$ are given by
\begin{eqnarray*}
\widetilde d_k&=&c_k, \quad  \quad k=0,\ldots,p+q,\\
\widetilde d_k&=&\sum_{i=0}^p \frac{\widetilde a_i}{\widetilde p_i^{k+1}}-\sum_{i=0}^q \widetilde b_i \sum_{j=0}^{k-1} \frac{\widetilde d_j}{\widetilde z_i^{k-j+1}},  \quad \quad k=p+q+1, p+q+2, \ldots
\end{eqnarray*}

\vskip 2mm

In order to improve the order of approximation, it is possible, in theory, to choose the $\widetilde a_i$'s, $\widetilde b_i$'s, $\widetilde p_i$'s and $\widetilde z_i$'s such that they satisfy the system of nonlinear equations
$$c_k=\sum_{i=0}^p \frac{\widetilde a_i}{\widetilde p_i^{k+1}}-\sum_{i=0}^q \widetilde b_i \sum_{j=0}^{k-1} \frac{\widetilde d_j}{\widetilde z_i^{k-j+1}}, \quad k=0, \ldots,2(p+q+1),$$
in which case we have
$$\widetilde R_{p,q}(t)-f(t)={\cal O}(t^{2(p+q+1)+1}).$$
Since $\widetilde R_{p,q}$ has a numerator and a denominator both of degree $p+q+1$, it will be identical to the Pad\'e
approximant $[p+q+1/p+q+1]_f$. Obviously, in practice, the solution of this system is not an easy task.

\subsubsection{Form 2}

 For the form \eqref{f2}, the condition \eqref{cond} is
$$\sum_{i=0}^p \frac{\widehat a_i}{1-\widehat p_it}=(c_0+c_1t+c_2t^2+\cdots)\sum_{i=0}^q\frac{\widehat b_i}{1-\widehat z_it}.$$
But, $1/(1-\widehat p_it)=1+\widehat p_it+\widehat p_i^2t^2+\cdots$, and a similar expansion for $1/(1-\widehat z_it)$. Thus, after replacement, the preceding relation becomes
$$\sum_{i=0}^p \widehat a_i (1+\widehat p_it+\widehat p_i^2t^2+\cdots)=(c_0+c_1t+c_2t^2+\cdots)
\sum_{i=0}^q \widehat b_i(1+\widehat z_it+\widehat z_i^2t^2+\cdots).$$

Identifying the coefficients of the identical powers of $t$ in both  sides, and taking into account the normalization condition leads to the system of equations for determining the coefficients $\widehat a_i$ and $\widehat b_i$
\begin{equation}
\label{coef2}
\left.
\begin{array}{l}
\displaystyle \sum_{i=0}^q \widehat b_i=1 \\
\displaystyle \sum_{i=0}^p \widehat a_i \widehat p_i^{k}-\sum_{i=0}^q \widehat b_i \sum_{j=0}^{k-1} c_j \widehat z_i^{k-j}=c_k, \qquad k=0,\ldots,p+q.
\end{array}
\right\}
\end{equation}

Again, the sum on $j$ is empty for $k=0$.

\vskip 2mm

The coefficients of the series expansion of $\widehat R_{p,q}(t)=\widehat d_0+\widehat d_1t+\widehat d_2t^2+\cdots$ are given by
\begin{eqnarray*}
\widehat d_k&=&c_k,  \quad \quad k=0,\ldots,p+q,\\
\widehat d_k&=&\sum_{i=0}^p \widehat a_i \widehat p_i^{k}-\sum_{i=0}^q \widehat b_i \sum_{j=0}^{k-1} \widehat d_j \widehat z_i^{k-j},  \quad \quad k=p+q+1, p+q+2, \ldots
\end{eqnarray*}

The order of approximation of $\widehat R_{p,q}$ can be improved as explained for $\widetilde R_{p,q}$.

\subsection{Barycentric Pad\'e--type approximants ({\sc bpta})}

Consider again the rational functions \eqref{f1} and \eqref{f2} and assume now that the coefficients $b_i$ in their respective denominators are arbitrarily chosen. Then, the coefficients of their numerators can be computed by solving the system \eqref{coef3}

\begin{equation}
\label{coef3}
\displaystyle \sum_{i=0}^p \frac{\widetilde a_i}{\widetilde p_i^{k+1}}=\sum_{i=0}^q \widetilde b_i \sum_{j=0}^{k-1} \frac{c_j}{\widetilde z_i^{k-j+1}}, \qquad k=0,\ldots,p
\end{equation}
for \eqref{f1}, or the system \eqref{coef4} for \eqref{f2}

\begin{equation}
\label{coef4}
\displaystyle \sum_{i=0}^p \widehat a_i \widehat p_i^{k}=\sum_{i=0}^q \widehat b_i \sum_{j=0}^{k-1} c_j \widehat z_i^{k-j}, \qquad k=0,\ldots,p.
\end{equation}

In both cases, the rational function $R_{p,q}$ which is obtained satisfies
$$R_{p,q}(t)-f(t)={\cal O}(t^{p+1}),$$
and, thanks to this property, it is called a {\it barycentric Pad\'e--type approximant}  (see \cite{pta}) and denoted {\sc bpta}.

\vskip 2mm

Similarly, a part of the numerator and/or a part of the denominator can be fixed thus leading to a {\it barycentric partial Pad\'e--type approximant} in the style of \cite{ppa}.

\vskip 2mm

The $b_i$'s could be chosen so that $R_{p,q}$ satisfies some additional properties as explained above.

\section{Partial fraction form ({\sc pfpa})}

Let us consider now the rational function
\begin{equation}
\label{pff}
R_{k,k+1}(t)=\sum_{i=0}^k \frac{a_i}{1-p_it}.
\end{equation}
It has a denominator of degree $k+1$ and a numerator of degree $k$. We want to compute the $a_i$'s and the $p_i$'s such that this rational function be identical to the Pad\'e approximant $[k/k+1]_f$ of the series $f$.
Such an approximant will be called a {\it partial fraction Pad\'e approximant} and denoted {\sc pfpa}.
It can be obtained by a slight variation of a method due to the French mathematician and hydraulics engineer
Gaspard Clair Fran\c{c}ois Marie Riche, Baron de Prony (Chamelet, 22 July 1755 - Asni\`eres--sur--Seine, 29 July 1839)
for interpolation by a sum of exponential functions \cite{pron}. This method is used in signal analysis and recovery (see, for example, \cite{pl1,pl2}). Applied to our case, this variant is as follows (see, for example, \cite[pp. 141--142]{cont}).

We want to have
$$\sum_{i=0}^k \frac{a_i}{1-p_it}=\sum_{i=0}^k a_i(1+p_it+p_i^2t^2+\cdots)=c_0+c_1t+c_2t^2+\cdots,$$
which leads, by identification of the powers of $t$ on both sides, to
\begin{equation}
\label{eqpr}
\sum_{i=0}^ka_ip_i^j=c_j, \quad j=0,\ldots,2k+1.
\end{equation}
The denominator of the rational function \eqref{pff} is
$$Q(t)=(1-p_0t)\cdots(1-p_kt)=b_0+b_1t+\cdots+b_{k+1}t^{k+1},$$
where $b_0=1$.
Let us first compute its coefficients. Multiply the first equation in \eqref{eqpr} (that is the equation for $j=0$) by $b_0$, the second one (that is corresponding to $j=1$) by $b_1$, and so on up to the $(k+2)$th equation (that is the equation for
$j=k+1$) by $b_{k+1}$, and sum them up. Begin again the same process starting from the second equation in \eqref{eqpr} (that is for $j=1$) which is multiplied by $b_0$, multiply the third equation by $b_1$, and so on up to the $(k+3)$th equation (that is the equation for $j=k+2$) by $b_{k+1}$, and sum them up. Continue the process until all equations in \eqref{eqpr}
have been used. We finally obtain
$$\sum_{j=0}^{k+1}b_j c_{j+n}=\sum_{j=0}^{k+1}b_j \sum_{i=0}^k a_i p_i^{j+n}, \quad n=0,\ldots,k,$$
which can be written as
$$\sum_{j=0}^{k+1}b_j c_{j+n}=\sum_{i=0}^k a_i p_i^n \sum_{j=0}^{k+1} b_j p_i^j=\sum_{i=0}^k p_i^n Q(p_i)=0, \quad n=0,\ldots,k.$$
Since $b_0=1$, the other coefficients $b_1,\ldots,b_{k+1}$ of the polynomial $Q$ are solution of the linear system
\begin{equation}
\label{fop}
\sum_{j=1}^{k+1} b_j c_{j+n}=-c_n, \quad n=0,\ldots,k.
\end{equation}
After solving this system, the zeros $p_0,\ldots,p_k$ of $Q$ can be computed (for example, by the $QR$ algorithm as the eigenvalues of the companion matrix of the coefficients of the polynomial $Q$) and, finally, $a_0,\ldots,a_k$ are obtained by solving the linear system consisting in the first $k+1$ equations in \eqref{eqpr}. Notice that this system is singular if the $p_i$'s are not distinct and, in this case, $[k/k+1]_f$ cannot be written under the form \eqref{pff} but, possibly, under
a partial fraction form involving powers in its denominator. This case arises if the denominator of $[k/k+1]_f$ has multiple zeros. It is easy to see that the system \eqref{fop} is identical to the system \eqref{bi} when $p=k$ and $q=k+1$, and after reversing the numbering of its coefficients.
There exist several procedures for improving the numerical stability of Prony's method \cite{potts0}, and
it is even possible to avoid the computation of the coefficients of $Q$ \cite{potts}.

\vskip 2mm

Let $c$ be the linear functional on the vector space of polynomials defined by
$$c(x^i)=c_i, \quad i=0,1,\ldots$$
Then, the system \eqref{fop} can be written
$$c(x^nQ(x))=0, \quad n=0,\ldots,k.$$
Thus, $Q$ is the polynomial of degree $k+1$ belonging to the family of {\it formal orthogonal polynomials} with respect to $c$. Such polynomials, introduced by Wynn \cite{pw}, play an important role in the algebraic and in the analytic theory of Pad\'e approximation \cite{pta,birk,hsvt}. Via formal orthogonal polynomials, Pad\'e approximants are also related to formal Gaussian quadrature procedures \cite{quad}. An interesting reference on the connection between these topics is
\cite{weiss}.

\vskip 2mm

Obviously, if the $p_i$'s all distincts, then $a_i=-p_iP(1/p_i)/Q'(1/p_i)$ where $P$ is the numerator of  $[k/k+1]_f$.

\vskip 2mm

If the $p_i$'s are arbitrary distinct points and if the $a_i$'s are solution of the first $k+1$ equations in \eqref{eqpr}, then
$R_{k,k+1}$ is only a Pad\'e--type approximant of $f$. In this case, the $p_i$'s could be chosen so that $R_{k,k+1}$ satisfies some additional properties.

\vskip 2mm

Partial fraction Pad\'{e} approximants can also be written as $\sum_{i=0}^k a_i/(p_i-t)$, and treated similarly.
Analogous forms can be derived for Partial Pad\'e approximants \cite{ppa}.

\begin{figure}[p]
\begin{center}
\includegraphics[width=0.5\textwidth]{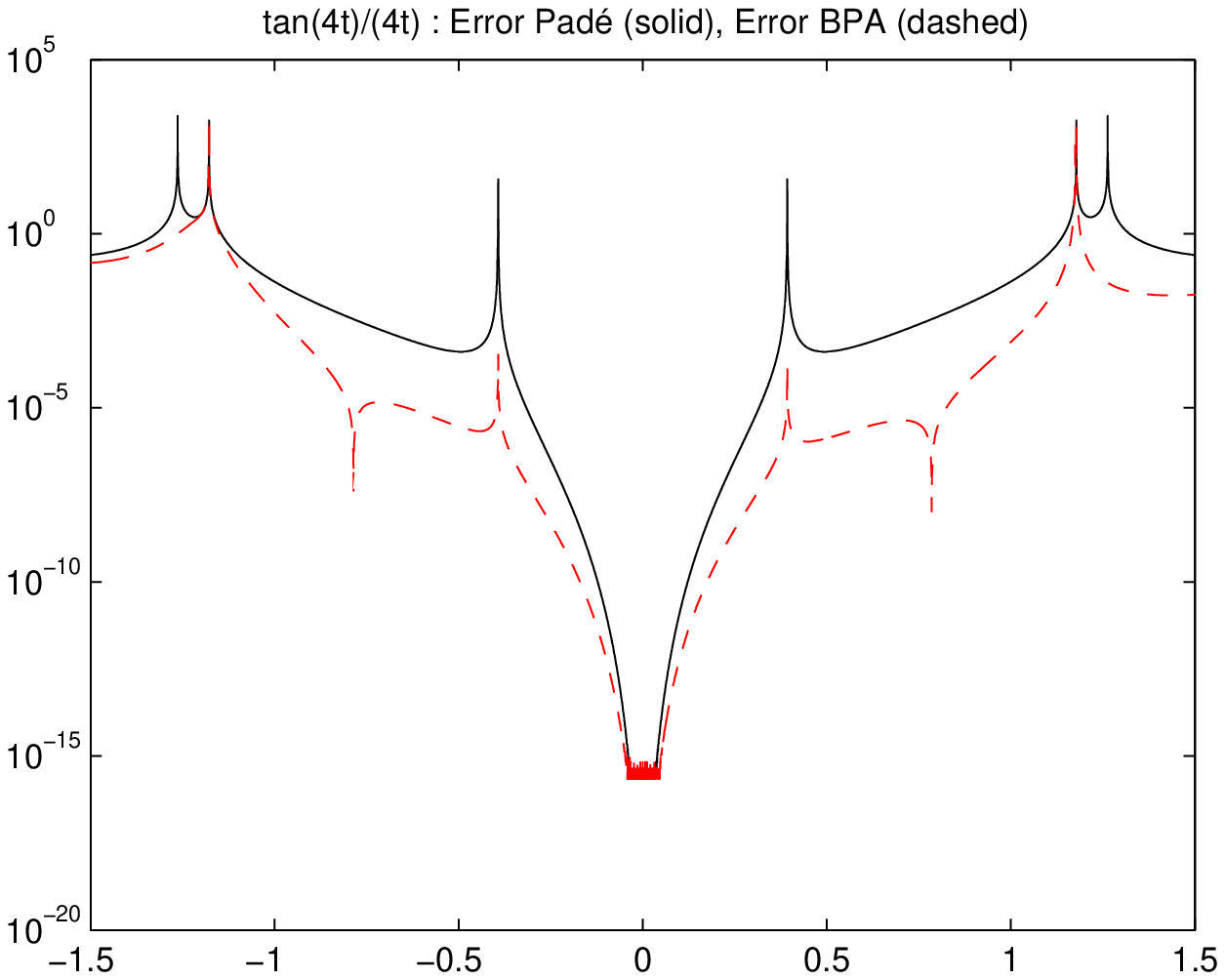}
\caption{Error of Pad\'{e} (solid) and barycentric Pad\'{e} (dashed) approximants for $\tan(4t)/(4t)$.}
\label{ftan1}
\vspace{0.5cm}
\includegraphics[width=0.5\textwidth]{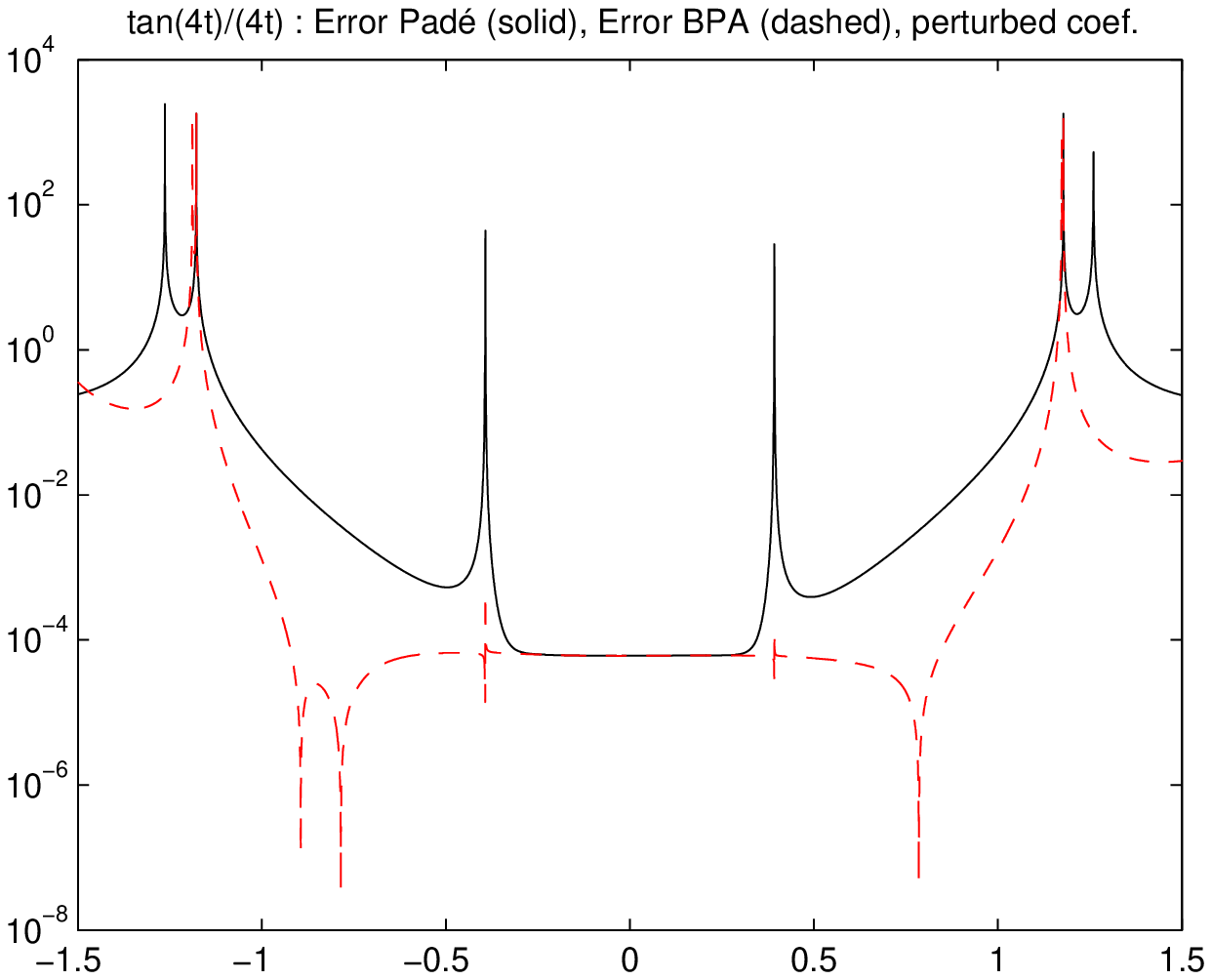}
\caption{Error of Pad\'{e} (solid) and barycentric Pad\'{e} (dashed) approximants for $\tan(4t)/(4t)$ with perturbed coefficients.}
\label{ftan2}
\end{center}
\end{figure}

\subsection{Numerical examples}

Let us now give some numerical examples showing the interest of the barycentric forms.
Both forms give similar results as expected.

\subsubsection{Example 1}

 We consider the following function, and its series expansion
$$f(t)=\frac{\tan(\omega t)}{\omega t}=1+\frac{1}{3}\omega^2t^2+\frac{2}{15}\omega^4t^4
+\frac{17}{315}\omega^6t^6+\frac{62}{2835}\omega^8t^8+\cdots$$
This function has poles at odd multiples of $\pi/(2\omega)$, and zeros at odd multiples of $\pi/\omega$, except at 0.

With $\omega=4$, $p=q=4$ and taking for the $\widetilde p_i$'s $\pm \pi/(2\omega), \pm 3\pi/(2\omega), 5 \pi/(2\omega)$ and the for $\widetilde z_i$'s $\pm \pi/\omega, \pm 3\pi/\omega, 5 \pi/\omega$, which are the five first poles and zeros of $f$ respectively, we obtain the results of Figure \ref{ftan1}.

\vskip 0.2cm
Adding a uniformly distributed random perturbation between $[-0.0001,+0.0001]$ to the $c_i$'s leads to the results of
Figure \ref{ftan2}. For $t \in [-1.5,+1.5]$, the error of the true Pad\'e approximant is in the interval $[3.5528 \times 10^{-5},2.8663 \times 10^{4}]$, and for the barycentric Pad\'e approximant computed either by the system \eqref{coef1} (form 1) of by the system \eqref{coef2} (form 2) it belongs to
$[4.8921 \times 10^{-8},1.8136 \times 10^{3}]$.

\begin{figure}[htb]
\begin{center}
\includegraphics[width=0.5\textwidth]{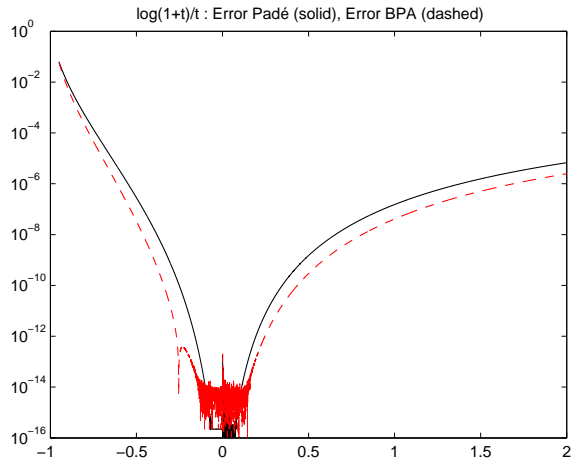}
\caption{Error of Pad\'{e} (solid) and barycentric Pad\'{e} (dashed) approximants for $\log(1+t)/t$.}
\label{flog1}
\vspace{-0.6cm}
\end{center}
\end{figure}
\subsubsection{Example 2}

We consider the series
$$f(t)=\frac{\log(1+t)}{t}=1-\frac{t}{2}+\frac{t^2}{3}-\frac{t^3}{4}+\cdots$$
which converges in the unit disk and on the unit circle except at the point $-1$ since there is a cut from
$-1$ to $-\infty$.
For $p=q=4$, the $p_i$'s equidistant in $[-10,-1]$ and the $z_i$'s equidistant in $[-10,-2]$, we obtain the results of Figure \ref{flog1}.
The numerical results highly depend on these choices.

\section{Conclusion}
This paper is an addition to the vast literature on Pad\'{e} approximation, and is only an introduction to these barycentric and partial fraction forms in order to show how to compute their coefficients. Their main features are discussed and some numerical experiments show the interest of these new representations. However, the main problem which remains to be studied is the influence of the choice of the free parameters involved in their construction, a choice related to the important issues of their robustness \cite{tref} and their stability \cite{bbam}.


\end{document}